\documentclass[12pt,a4paper]{amsart}
\usepackage{amssymb}
\usepackage{multirow}
\usepackage{graphicx}
\usepackage{amsmath,amsbsy,amsfonts,amsthm,latexsym,amsopn,amstext,amsxtra,euscript,amscd,mathrsfs}
\usepackage{cite}

\theoremstyle{plain}

\theoremstyle{definition}

\theoremstyle{remark}

\numberwithin{equation}{section}

\numberwithin{table}{section}

\numberwithin{figure}{section}

\begin{document}

\title[An integral approach to the Gardner-Fisher and untwisted Dowker sums]{An
integral approach to the Gardner-Fisher and untwisted Dowker sums}


\author{Carlos M. da Fonseca}
\address{Department of Mathematics, Kuwait University, Safat 13060, Kuwait}
\email{carlos@sci.kuniv.edu.kw}

\author{M. Lawrence Glasser}
\address{Department of Physics, Clarkson University, Potsdam, NY 13699-5820, USA}
\email{laryg@clarkson.edu}

\author{Victor Kowalenko}
\address{ARC Centre of Excellence for Mathematics and Statistics of Complex Systems, Department of Mathematics and Statistics,
The University of Melbourne, Victoria 3010, Australia}
\email{vkowa@unimelb.edu.au}

\subjclass[2000]{30B10,11M06,33E20,33B10}

\date{March 1, 2015}

\keywords{digamma function, Gardner-Fisher sum, generalized cosecant numbers, symmetric polynomial, trigonometric power sum,
untwisted Dowker sum, zeta function}

\date{\today}

\begin{abstract}
We present a new and elegant integral approach to computing the Gardner-Fisher
trigonometric power sum, which is given by
$$
S_{m,v}=\left(\frac
\pi{2m}\right)^{2v}\sum_{k=1}^{m-1}\cos^{-2v}\left(\frac{k\pi}{2m}\right)\,
,
$$
where $m$ and $v$ are positive integers. This method not only
confirms the results obtained earlier by an empirical method, but it
is also much more expedient from a computational point of view. By
comparing the formulas from both methods, we derive several new
interesting number theoretic results involving symmetric polynomials
over the set of quadratic powers up to $(v-1)^2$ and the generalized
cosecant numbers. The method is then extended to other related
trigonometric power sums including the untwisted Dowker sum.
By comparing both forms for this important sum, we derive new
formulas for specific values of the N\"{o}rlund polynomials.
Finally, by using the results appearing in the tables, we consider
more advanced sums involving the product of powers of cotangent and
tangent with powers of cosecant and secant respectively.

\end{abstract}

\maketitle

\section{Introduction}

In 1969 Gardner \cite{Ga1969} stated that the finite sum of inverse
powers of cosines
\begin{equation}\label{one}
S_{m,v}=\left(\frac{\pi}{2m}\right)^{2v}\sum_{k=1}^{m-1}\cos^{-2v}\left(\frac{k\pi}{2m}\right) \,
,
\end{equation}
where $m$ and $v$ are positive integers, emerges during the calculation of the $v$th cumulant of a
certain quadratic form in $m$ independent standardized normal variates. Although he was able to show that
$$
\lim_{m\rightarrow\infty}S_{m,v}=\zeta(2v)\, ,
$$ where $\zeta$ is the Riemann zeta function, he posed the problem of whether it was possible to obtain a
``simpler" closed form expression for $S_{m,v}$, for all $m$ and $v$.

In solving this problem, Fisher \cite{Fi1971,Kl1990} observed that the above sum could also be written as
\begin{equation}\label{onea}
S_{m,v}=\left(\frac
\pi{2m}\right)^{2v}\sum_{k=1}^{m-1}\sin^{-2v}\left(\frac{k\pi}{2m}\right)\, ,
\end{equation}
and then devised an ingenious generating function approach for the more general trigonometric power series given by
$$
Q_{m,v}(\delta)=\sum_{k=1}^{m-1}\sin^{-2v}\left(\frac{k\pi+\delta}{2m}\right)\,
.
$$
Consequently, he found that
$$
S_{m,1}=\frac{\pi^2}6\left(1-\frac 1{m^2}\right) \quad
\mbox{and}\quad S_{m,2}=\frac{\pi^4}{90}\left(1+\frac 5{2m^2}-\frac
7{2m^4}\right)\, ,
$$
while for large $m$, he obtained
$$
S_{m,v}=\zeta(2v)+\frac v{12m^2}\, \zeta(2v-2)+O\left(m^{-4}\right)\, .
$$
Although Gardner first obtained $S_{m,v}$, Fisher's solution is largely responsible for the continuing interest in the sum.
Consequently, the sum will be henceforth referred to as the Gardner-Fisher sum.

\section{Generalized cosecant numbers}
In a recent work one of us \cite{K2011b} showed that the Gardner-Fisher sum can be
expressed as
$$
S_{m,v}=\sum_{j=0}^\infty
c_{2v,j}\left(\frac\pi{2m}\right)^{2j}\sum_{k=1}^{m-1} k^{2j-2v}\, ,
$$
where the coefficients $c_{\rho,k}$ are a generalization of the
cosecant numbers $c_{k}$ in Ref.\ \cite{K2011a} and are given by
\begin{equation}\label{kow1}
c_{\rho,k}= (-1)^k \sum_{\substack{n_1,n_2,n_3,\dots,n_k
\\\sum_{i=1}^k i\, n_i=k}}^{k,\lfloor k/2 \rfloor, \lfloor k/3
\rfloor,\dots,1} (-1)^N (\rho)_N \prod_{i=1}^k \Bigl(
\frac{1}{(2i+1)!} \Bigr)^{n_i}  \frac{1}{n_i!}\, .
\end{equation}

The above formula has been derived by applying the partition method
for a power series expansion to $x^{\rho}/\sin^{\rho}x$. That is,
the generalized cosecant numbers are the coefficients of the power
series for this function. The partition method for a power series
expansion was first introduced in Ref.\  \cite{K1994}, but more recently
it has been developed further in Refs.\ \cite{K2011a}-\cite{K2009} and
\cite{K2012}. To calculate the generalized cosecant number $c_{\rho,
k}$ via \eqref{kow1}, we need to determine the specific contribution
made by each integer partition that sums to $k$. For example, if we
wish to determine $c_{\rho,5}$ we require all the contributions made
by the seven partitions that sum to 5 or those listed in the first
column of Table\ \ref{table1}. Depending on the function that is
being studied, each element or part in a partition is assigned a
specific value. In the case of $x^{\rho}/\sin^{\rho}x$, each element
$i$ is assigned a value of $(-1)^{i+1}/(2i+1)!$. Moreover, if an
element occurs $n_i$ times in a partition or possesses a
multiplicity of $n_i$, then we need to take $n_i$ values or
$(-1)^{(i+1)n_i}/((2i+1)!)^{n_i}$. Table\ \ref{table1} displays the
multiplicities of all elements in the partitions that sum to 5.

\begin{table}
\small
\begin{tabular}{|c|c c c c c|} \hline
Partition & $n_1$ & $n_2$ & $n_3$ & $n_4$ & $n_5$    \\ \hline
$\{5\}$ & & & & & $ 1 $  \\
$\{4,1\}$ & 1 &  & & 1 &   \\
$\{3,2\}$ & & 1 & 1 & &    \\
$\{3,1,1\}$ & 2 & & 1 & &  \\
$\{2,2,1\}$ & 1 & 2 & & &   \\
$\{2,1,1,1\}$ & 3 & 1& & &   \\
$\{1,1,1,1,1\}$ & 5 & & & & \\
\hline
\end{tabular}
\vspace{0.3cm}
\caption{Multiplicities of the partitions whose elements or parts sum to 5}
\label{table1}\end{table}

Associated with each partition is a multinomial factor that is
determined by taking the factorial of the total number of elements
in the partition $N(=\sum_{i}n_i)$ and dividing by the factorials of
all the multiplicities. For the partition $\{2,1,1,1\}$ in Table\
\ref{table1}, we have $n_1=3$ and $n_2=1$, while all other
multiplicities do not contibute or vanish. Hence, the multinomial
factor for this partition is $4!/(3!\, 1!)=4$. When the function is
accompanied by a power, another modification must be made. Each
partition must be multiplied by Pochhammer factor of
$\Gamma(N+\rho)/\Gamma(\rho) N!$ or $(\rho)_N/N!$. For $\rho=1$,
this simply yields unity and thus, the multinomial factor is
unaffected. Then the generalized cosecant numbers reduce to the
cosecant numbers $c_k$ \cite{K2011a}, which are defined as
$$
c_k = \frac{(-1)^{k+1}}{(2k)!}(2^{2k}-2) B_{2k}=2
\left(1-2^{1-2k}\right) \frac{\zeta(2k)}{\pi^{2k}} \, ,
$$
and $B_{2k}$ represent the Bernoulli numbers.

In \eqref{kow1} the product is concerned with the evaluating the
contribution made by each partition based on the values of the
multiplicities, while the sum is concerned with all partitions.
Thus, the sum covers the range of values for the multiplicities. For
example, $n_1$ attains a maximum value of $k$, which corresponds the
partition with $k$ ones, while $n_2$ attains a maximum value of
$[k/2]$, which corresponds to the partition with $[k/2]$ twos in it.
Here, $[k/n]$ denotes the floor function or the largest integer less than
or equal to $k/n$. Therefore for odd values, the partition with $[k/2]$ twos
would also have $n_1=1$. Hence, we see that $n_i$ can only attain a maximum value of
$[k/i]$, which becomes the upper limit for each multiplicity in \eqref{kow1}.
Furthermore, a valid partition must satisfy the constraint given by
$\sum_{i=1}^k in_i=k$.

As an example, let us consider the evaluation of $c_{\rho,5}$. According to Table\ \ref{table1} there are
seven partitions whose contributions must be evaluated. By following the steps given above, we find that
in the order in which the partitions appear in Table\ \ref{table1}, \eqref{kow1} yields
\begin{eqnarray}\label{kow2}
c_{\rho,5}&= (\rho)_1 \;\frac{1}{11!} - \frac{(\rho)_2}{2!}\;
\frac{2!}{1!  1!}\; \frac{1}{3! \cdot 9!} -\frac{(\rho)_2}{2!} \;
\frac{2!}{1! 1!}\; \frac{1}{5!  7!} + \frac{(\rho)_3}{3!}\;
\frac{3!}{1! 2!} \; \frac{1}{3!^2  7!}
\nonumber\\
& +\;\; \frac{(\rho)_3}{3!}\; \frac{3!}{1! 2!}\; \frac{1}{5!^2 3!} -
\frac{(\rho)_4}{ 4!}\, \frac{4!}{1! 3!} \; \frac{1}{3!^3 5!} +
\frac{(\rho)_5}{5!}\; \frac{5!}{5!} \;\frac{1}{3!^5} \;\;.
\end{eqnarray}
From the above result we see that $c_{\rho,5}$ is a fifth order
polynomial in $\rho$. In fact, the highest order term comes from the
partition with $k$ ones in it, which produces the term with
$(\rho)_k$. Hence, we observe that the generalized cosecant numbers
are polynomials of order $k$. Another interesting property is that
all the contributions from partitions with the same number of
elements in them possess the same sign, which toggles according to
whether there is an even or odd number of elements in the
partitions. Furthermore, \eqref{kow2} can be simplified further,
whereupon one arrives at the result appearing in the sixth row below
the headings in Table \ref{table2}.

\begin{table}
\begin{center}
\small
\begin{tabular}{|c|l|} \hline
$k$ & $c_{\rho,k}$    \\ \hline
$0$ &  $ 1 $  \\ [0.1 cm]
$1$   &  $\frac{1}{3!} \; \rho $ \\ [0.1 cm]
$2$  &  $ \frac{2}{6!} \; \bigl(2\rho+5\rho^2 \bigr)$  \\ [0.1 cm]
$3$ &  $ \frac{8} {9!}\;\bigl( 16 \rho+42 \rho^2+35 \rho^3 \bigr) $ \\ [0.1 cm]
$4$ &  $ \frac{2}{3 \cdot 10!} \; \bigl( 144 \rho+404\rho^2+420 \rho^3+175 \rho^4 \bigr)$ \\ [0.1 cm]
$5$ & $ \frac{4}{3 \cdot 12!} \; \bigl( 768 \rho+2288 \rho^2+2684 \rho^3+1540 \rho^4+ 385 \rho^5 \bigr)$ \\ [0.1 cm]
$6$ & $ \frac{2}{9 \cdot 15!} \; \bigl( 1061376 \rho+3327594 \rho^2+4252248 \rho^3+2862860 \rho^4+1051050
\rho^5  +175175 \rho^6 \bigr)$ \\ [0.1 cm]
$7$  & $ \frac{1} {27 \cdot 15!} \; \bigl( 552960 \rho+ 1810176 \rho^2+ 2471456 \rho^3 + 1849848 \rho^4+ 820820 \rho^5 + 210210 \rho^6$ \\
& $ +25025 \rho^7 \bigr)$ \\ [0.1 cm]
$8$ & $ \frac{2}{45 \cdot 18!} \bigl( 200005632 \rho + 679395072 \rho^2+ 978649472 \rho^3 + 792548432 \rho^4 + 397517120 \rho^5 $ \\
& $+ 125925800 \rho^6 + 23823800 \rho^7+ 2127125 \rho^8 \bigr)$ \\ [0.1 cm]
$9$ & $ \frac{4}{81 \cdot 21!} \bigl(129369047040 \rho+ 453757851648 \rho^2 + 683526873856 \rho^3 +  589153364352 \rho^4 $ \\
& $ + 323159810064 \rho^5 + 117327450240 \rho^6 + 27973905960 \rho^7 + 4073869800 \rho^8 +  282907625 \rho^9 \bigr)$ \\[0.1 cm]
$10$ & $\frac{2}{6075 \cdot 22!} \bigl( 38930128699392 \rho + 140441050828800 \rho^2 + 219792161825280 \rho^3 + 199416835425280 \rho^4 $ \\
& $ + 117302530691808 \rho^5 + 47005085727600 \rho^6 + 12995644662000 \rho^7 + 2422012593000 \rho^8$ \\
& $ + 280078548750 \rho^9 + 15559919375 \rho^{10} \bigr)$  \\[0.1 cm]
$11$ & $ \frac{8}{243 \cdot 25!} \bigl(494848416153600 \rho + 1830317979303936 \rho^2 + 2961137042841600 \rho^3 $ \\
& $ + 2805729689044480 \rho^4 + 1747214980192000 \rho^5 + 755817391389984 \rho^6 + 232489541684400 \rho^7$ \\
& $ + 50749166067600 \rho^8 + 7607466867000 \rho^9 + 715756291250 \rho^{10} + 32534376875 \rho^{11} \bigr)$ \\[0.1 cm]
$12$ & $ \frac{2}{2835 \cdot 27!} \bigl( 1505662706987827200 \rho + 5695207005856038912 \rho^2 + 9487372599204065280 \rho^3$ \\
& $ + 9332354263294766080 \rho^4 + 6096633539052376320 \rho^5 + 2806128331871953088 \rho^6 $ \\
& $ + 937291839756592320 \rho^7 + 229239926321406000 \rho^8 + 40598842049766000 \rho^9 $ \\
& $ + 5005999501002500 \rho^{10} + 390802935022500 \rho^{11} + 14803141478125 \rho^{12} \bigr) $ \\[0.1 cm]
$13$ & $ \frac{232}{81 \cdot 30!} \bigl( 844922884529848320 \rho + 3261358271400247296 \rho^2 +5576528334428209152 \rho^3$ \\
& $+5668465199488266240 \rho^4 + 3858582205451484160 \rho^5 + 1870620248833400064 \rho^6 $ \\
& $ +667822651436228288 \rho^7 + 178292330746770240 \rho^8 + 35600276746834800 \rho^9 $ \\
& $ 5225593531158000 \rho^{10} + 539680243602500 \rho^{11} + 35527539547500 \rho^{12} + 1138703190625 \rho^{13} \bigr) $ \\[0.1 cm]
$14$ & $ \frac{2}{1215 \cdot 30!}\; \bigl(138319015041155727360 \rho + 543855095595477762048 \rho^2 + 952027796641042464768 \rho^3 $ \\
& $ +996352286992030556160 \rho^4 + 703040965960031795200 \rho^5 + 356312537387839432192 \rho^6 $ \\
& $ +134466795172062184832 \rho^7 + 38526945410311117760 \rho^8 + 8436987713444690400 \rho^9 $\\
& $+ 1404048942958662000 \rho^{10} + 173777038440005000 \rho^{11} + 15258232341852500 \rho^{12} $\\
& $ + 858582205731250 \rho^{13} +23587423234375 \rho^{14} \bigr)$ \\ [0.1 cm]
$15$ & $\frac{1088}{729 \cdot 35!} \; \bigl( 562009739464769840087040 \rho + 2247511941596311764074496 \rho^2 $ \\
& $4019108379306905439830016 \rho^3 + 4317745925208072594259968 \rho^4 $ \\
& $ + 3145163776677939429416960 \rho^5 +1656917203539032341530624 \rho^6 $ \\
& $ + 655643919364420586023424 \rho^7 + 199227919419039256217472 \rho^8 $ \\
& $ + 46995751664475880185920 \rho^9 + 8614026107092938211680 \rho^{10} +1214778349162323946000 \rho^{11} $ \\
& $+ 128587452922193265000 \rho^{12} + 9720180867524627500 \rho^{13} + 472946705787806250 \rho^{14}$ \\
& $ +11260635852090625 \rho^{15} \bigr) $ \\[0.1 cm]
\hline
\end{tabular}
\end{center}
\normalsize
\vspace{0.3cm}
\caption{Generalised cosecant numbers $c_{\rho,k}$ up to $k=15$}
\label{table2}\end{table}

Table\ \ref{table2} displays the first fifteen generalized cosecant
numbers obtained by determining the multiplicities of all the
partitions that sum to each order $k$ and then evaluating the sum of
their contributions according to steps given above. Beyond
$k\!=\!10$, the partition method for a power series expansion
becomes laborious due to the exponential increase in the number of
partitions. To circumvent this problem, a general computing
methodology is required to determine higher order coefficients via
the partition method. This methodology, which is based on
representing all the partitions that sum to a specific order as a
tree diagram and invoking the bivariate recursive central partition
(BRCP) algorithm, is presented in Ref.\  \cite{K2012}. The general
expressions for the coefficients can then imported into a
mathematical software package such as Mathematica \cite{W1992}
whereupon its symbolic routines yield the final values presented in
Table \ref{table2}.

For the special case, where $\rho$ is an even integer (the case of interest here),
the generalized cosecant numbers satisfy the following recurrence relation:
\begin{equation}\label{kow2a}
c_{2n+2,k+1} = \frac{(2k+2-2n)}{2n} \, \frac{(2k +1-2n)}{(2n+1)}\, c_{2n, k+1} + \frac{2n}{2n+1}\, c_{2n,k}\,.
\end{equation}
This equation is obtained by introducing the power series expansion
for $x^{2n}/\sin^{2n}x$ into (27) of Ref.\ \cite{K2011b}, which is
\begin{equation} \label{kow2b}
\frac{d^2}{dx^2} \, \frac{1}{\sin^{2n} x}= \frac{2n}{\sin^{2n}x} +
\frac{2n(2n+1) \cos^2 x}{\sin^{2n+2} x} \;.
\end{equation}
Then one equates like powers of $x$. For the special case of $n=1$, the numbers are
related to the cosecant-squared numbers in Ref.\ \cite{K2011a}. Consequently, we find that
\begin{equation} \label{kow3}
c_{2,k}= 2 (2k-1) \, \frac{\zeta(2 k)}{\pi^{2 k}} \;.
\end{equation}
Furthermore, the $k$-dependence in \eqref{kow2a} can be removed by deducing that the
$c_{2n,k}$ can be expressed generally as
\begin{equation}\label{kow3a}
c_{2n,k} = 4\sum_{j=0}^{n-1} \frac{\Gamma(k-j)}{\Gamma(k-n+1)}\,
\frac{\Gamma(k-j+1/2)}{\Gamma(k-n+1/2)}\, \frac{1}{\Gamma(n)}\,
\frac{\Gamma(j+1/2)}{\Gamma(n+1/2)}\, C(n,j) \, \frac{\zeta(2k-2j)}{\pi^{2k-2j}}\;,
\end{equation}
where
\begin{equation}\label{kow3b}
C(n,j) = C(n-1,j) + \frac{(n-1)^2}{j-1/2}\, C(n-1,j-1) \;\;,
\end{equation}
and $C(1,0) \!=\! 2$ from equating \eqref{kow3a} with \eqref{kow3}. For $n \!=\! 1$, we find that
$C(n,1)= 2 \sum_{l=1}^{n-1} l^2 \!=\! 4 B_3(n)/3$, where $B_k(x)$ denotes a Bernoulli polynomial.
Introducing this result into \eqref{kow3b} yields $C(n,2)$, which is given by
\begin{equation}\label{kow3c}
C(n,2)=\frac{1}{135}\, n(n-1)(n-2)(2n-1)(2n-3)(5n+1)= \frac{1}{6} \,(2n-4)_4\, c_{2n,2} \,.
\end{equation}
Similarly, if one introduces \eqref{kow3c} into \eqref{kow3b}, then one obtains
\begin{equation}\label{kow3d}
C(n,3)=\frac{1}{60} \,(2n-6)_6\, c_{2n,3} \,.
\end{equation}
Therefore, we see that the $C(n,j)$ are related to the generalized cosecant numbers. In the next section
we shall also see that the generalized cosecant numbers are related to the symmetric polynomials over
the set of positive square integers.

\section{Integral approach}
In this section we present a totally different approach to
evaluating the Gardner-Fisher sum from others who have studied
similar trigonometric power sums, \cite{Be2002}-\cite{FK2013},
\cite{Ga2007}-\cite{GC1989}, \cite{G1968}, \cite{Kl1990} and
\cite{M2011b}- \cite{MT2012}.  The main advantage of the integral
approach presented here is that it yields a final form for the sum
that is more expedient from a computational point of view  than the
other references. In addition, except for Ref.\ \cite{K2011b}, it is both
more informative and compact than the other references. It also
provides an independent corroboration of the results in
Ref.\ \cite{K2011b}, which were obtained by empirical means. As a
consequence, a comparison of the results can be undertaken.

We begin this section by differentiating No.\ 8.365(10) in
Ref.\ \cite{GR1994}, which yields
\begin{equation}\label{three-one}
  \pi^2\csc^2(\pi x) = \psi'(x)+\psi'(1-x)\; .
\end{equation}
By introducing the series form for the derivative of the digamma
function, viz.\ $\psi'(x)=\sum_{n=0}^{\infty}1/(n+x)^2$ from No.\
8.363(8) of Ref.\ \cite{GR1994}, we replace the summand in
\eqref{three-one} by the integral representation for the gamma
function. On interchanging the order of the sum and the integral, we
obtain
$$
 \pi^2 \csc^2(\pi x) = \int_0^{\infty}\sum_{n=0}^{\infty}\left( e^{-(n+x)t}+e^{-(n+1-x)t}\right) t\, dt\;.
 $$
Next we evaluate the sums over $n$ via the geometric series. This yields
$$
 \pi^2 \csc^2(\pi x)  = \int_0^{\infty} \left( \frac{e^{-xt}}{1-e^{-t}}+\frac{e^{-(1-x)t}}{1-e^{-t}}\right) t\, dt \;.
$$
Now we make the change of variable $u=e^{-t}$, thereby arriving at
\begin{equation}
  \pi^2 \csc^2(\pi x) = -\int_0^1 \frac{\ln u}{1-u} \Bigl( u^x+u^{1-x}\Bigr) \frac{du}{u}. \label{eq1}
\end{equation}

With the aid of the following identity:
$$
\prod_{n=1}^{v-1}\left(\partial_z^2+4n^2\right)\csc^2z=(2v-1)!\csc^{2v}z\,
, \quad \mbox{for $v=1,2,3,\ldots$}\;,
$$
which is obtained by replacing the $\cos^2 x$ by $1 \!-\!\sin^2 x$ in \eqref{kow2b},
thereby yielding
$$
\Bigl( \frac{d^2}{dx^2} + 4 n^2 \Bigr) \csc^{2n}x = 2n(2n+1) \csc^{2n+2} x\;,
$$ and then multiplying this resulting equation by successive values of $n$, we find after an elementary
change of variable that \eqref{eq1} can be expressed as
$$
\csc^{2v}\left(\frac{k\pi}{2m}\right)=-\frac{4^vm^2}{(2v-1)!\pi^2}\int_0^1
\frac{\ln u}{1-u^{2m}} \prod_{n=1}^{v-1}\left(
\frac{m^2}{\pi^2}\ln^2 u
+n^2\right)\left(u^k+u^{2m-k}\right)\frac{du}{u}\, .
$$
Carrying out the summation over $k$ in the Gardner-Fisher sum or \eqref{onea} yields
\begin{equation}\label{eq2}
\sum_{k=1}^{m-1}\csc^{2v}\left(\frac{k\pi}{2m}\right)=-\frac{4^vm^2}{(2v-1)!\pi^2}\int_0^1\prod_{n=1}^{v-1}
\left(\frac{m^2}{\pi^2}\ln^2u+n^2\right)\frac{(1-u^{m-1})}{(1-u^m)(1-u)}\,
du\, .
\end{equation}

If we introduce Newton's identities for symmetric polynomials \cite{wi2014}, then \eqref{eq2} becomes
\begin{equation}\label{eq3}
S_{m,v}= -\frac{4^v}{(2v-1)!}
\sum_{n=0}^{v-1}s(v,n)\left(\frac{\pi}{m}\right)^{2n-2v}\int_0^1
\frac{(1-u^{m-1})}{(1-u^m)(1-u)}\, \ln^{2v-2n-1}u\, du ,
\end{equation}
where $s(v,n)$ represents the $n$th elementary symmetric polynomial obtained
by summing quadratic powers, viz.\ $1^2,2^2,\ldots,(v-1)^2$. That
is,
$$s(v,n)=\sum_{1 \leq i_1<i_2<\cdots<i_n<v-1} x_{i_1} x_{i_2} \cdots
x_{i_n},$$ where $x_{i_1} <x_{i_2}< \cdots <x_{i_n}$ and each
$x_{i_j}$ is equal to at least one value in the set
$\left\{1,2^2,3^2,\dots,(v-1)^2 \right\}$. In particular, for the
three lowest values of $n$, they are given by
$$ s(v,0)=1 \;,  \quad s(v,1)=(v-1)v(2v-1)/6 \; ,$$
and $$ s(v,2)= \frac{(5 v+1)}{4 \cdot 6!} \; (2v-4)_5 \; ,$$
while for the three highest values of $n$, they are given by
$$s(v,v-1) = (v-1)!^2 \;, \quad
s(v,v-2)=  (v-1)!^2\,\bigl( \zeta(2) -\zeta(2,v) \bigr) \;,$$
and $$ s(v,v-3) = \frac{(v-1)!^2}{2} \Bigl( (\zeta(2) -\zeta(2,v))^2+
\zeta(4,v)-\zeta(4) \Bigr) \; .$$
The integral $I$ over $u$ can now be evaluated by decomposing the
denominator of the integrand and applying No.\ 4.271(4) from
Ref.\ \cite{GR1994}. Then we find that  $$I =
\Gamma(2v-2n)\zeta(2v-2n) (m^{-2(v-n)}-1)\, .$$  Introducing this
result into \eqref{eq3}, we finally arrive at
\begin{equation}\label{three-three}
S_{m,v}=\frac{1}{(2v-1)!}\sum_{n=0}^{v-1} \Bigl(\frac{\pi}{m}\Bigr)^{2n} \, s(v,n)\Gamma(2v-2n) \,\zeta(2v-2n)
\left( 1-\frac{1}{m^{2v-2n}}\right)\, ,
\end{equation}

The above result for the Gardner-Fisher sum is computationally expedient because it can be
implemented as a one line instruction in Mathematica \cite{W1992}. With the aid of
the SymmetricPolynomial and the gamma and zeta function routines, \eqref{three-three} can be
expressed in Mathematica as
\vspace{0.3cm}

$ \rm{S[m_{-},v_{-}]:= (1/((2\;v-1)!))Sum[(Pi/m)^\wedge(2\;n)}$

$ \rm{SymmetricPolynomial[n,Table[k^\wedge 2,{k,1,v-1}]]} $

$\rm{Gamma[2\;v-2\;n] \; Zeta[2\;v-2\;n](1-1/(m^{\wedge}(2\;v-2\;n))),\{n,0,v-1\}]} \quad \quad.
$
\vspace{0.3cm}
\newline
If we type in the instruction \newline

F[m,n]:= Table[Simplify[S[m,v]],\{v,1,n\}] \hspace{0.5cm},
\vspace{0.3 cm}
\newline
then we can construct a table of values of the Gardner-Fisher sum. In fact, by typing \newline

Timing[Grid[Partition[F[m,15],1]]] \hspace{0.5 cm},
\vspace{0.3 cm}
\newline
one obtains the first 15 values of the Gardner-Fisher sum printed
out in table form. For a Venom Blackbook 17S Pro High Performance
laptop with 8 Gb RAM and equipped with Mathematica 10.1 these take
less than 0.16 seconds to appear on the screen. Table\ \ref{table3}
presents the output generated by the above instruction. We have not
only been able to present a much greater number of results than
Table\ 2 of  Ref.\ \cite{K2011b}, but we have also corrected the
typographical error occurring in the last term of the second result,
which should be $-7/2m^4$, not $-7/2m^2$. The results have been
presented in terms of the zeta function, thereby maintaining
consistency with Gardner's limit mentioned in the introduction.
Moreover, the results can be expressed as
$$
S_{m,v}= (2m^2)^{-v}\, (m^2-1)\, p_{v-1}(m^2) \;,
$$
where $p_{v}(x)$ is a polynomial of degree $v$ with the coefficients $p_{v,j}$ satisfying
\begin{equation}\label{three-four}
p_{v,v}=2^{v+1}\, \zeta(2v+2);\;, \quad p_{v,v-1}= 2^{v+1}
\zeta(2v+2) + 2^{v-1} \pi^2 \zeta(2v)/3 \;,
\end{equation}
$$
p_{v,j-1}=p_{v,j}+ \frac{2^{v+1}}{(2v+1)!}\,
\pi^{2v+2-2j}\,s(v+1,v+1-j) \Gamma(2j) \,\zeta(2j)\;,
$$
and
\begin{equation}\label{three-fourb}
p_{v,0}= \frac{2^{v+1}}{(2v+1)!}\, \sum_{j=0}^{v} \pi^{2j} \,
s(v+1,j) \Gamma(2v+2-2j)\,\zeta(2v+2-2j) \;.
\end{equation}

\begin{table}
\small
\begin{tabular}{|c|l|} \hline
$v$ & $(2m^2)^{v}S_{m,v}/ (m^2-1)$    \\ \hline
$1$ & $2\, \zeta(2)$ \\ [.1cm]
$2$  &  $2\, \zeta(4)  (7+2 m^2)$  \\ [.1cm]
$3$ &  $\zeta(6) (71+29 m^2+8 m^4)$ \\ [.1cm]
$4$ & $ \frac{2}{3} \, \zeta(8)(521+251 m^2+104 m^4+24 m^6)$ \\ [.1cm]
$5$ & $  \zeta(10) (1693+901 m^2+450 m^4+164 m^6+32 m^8)$ \\ [.1cm]
$6$ & $ \frac{1}{691}\, \zeta(12) (5710469+3253469 m^2+1815032 m^4+821182 m^6+262624 m^8+44224 m^{10})$ \\ [.1cm]
$7$ & $\frac{1}{30 }\, \zeta(14) (1212457+726457 m^2+436531 m^4+224891 m^6+91472 m^8+25952 m^{10}+3840 m^{12})$\\ [.1cm]
$8$ & $ \frac{2}{10851}\, \zeta(16) (1074010337+669172337 m^2+424359179 m^4+238674979 m^6+112425856 m^8$ \\
& $ +41352256 m^{10}+10528128 m^{12}+1388928 m^{14})$ \\ [.1cm]
$9$ & $ \frac{1}{219335}\, \zeta(18) (212920335247+136906047247 m^2+90463687339 m^4+54230674609 m^6$ \\
& $ +28237533526 m^8+12198655216 m^{10}+4088714368 m^{12}+943572928 m^{14}+112299520 m^{16})$ \\ [.1cm]
$10$ & $\frac{1}{3666831 } \, \zeta(20)(17471801743019+11530685023019 m^2+7871509097999 m^4 $\\
 & $ +4954255382249 m^6+2775847146014 m^8+1339874175764 m^{10}+533660885024 m^{12}$ \\
 & $ +164177309024 m^{14}+34637202944 m^{16}+3754834944 m^{18})$ \\ [.1cm]
$11$ &  $\frac{22}{51270780 }\, \zeta(22) (54543344015461+36797151215461 m^2 +25796160434461 m^4 $ \\
& $+16874095108441 m^6+10001178676906 m^8+5238677224486
    m^{10} +2352811821856 m^{12} $ \\
& $ +868592118736 m^{14}+246835895296 m^{16}+47953445376 m^{18} + 4772843520 m^{20})$   \\ [.1cm]
$12$ & $\frac{1}{3545461365 } \, \zeta(24)(407813841938063843+280376888294063843 m^2 $\\
 & $+200956537702208843 m^4+135642820502091743 m^6+84063888200850398
 m^8$ \\
 & $+46897459147087298 m^{10}+23026002351214568 m^{12}+9661082872248968 m^{14}$ \\
 & $+3322054926238208 m^{16}+876937263364608 m^{18}+157850700656640 m^{20}$ \\
& $ +14522209751040 m^{22})$ \\ [.1cm]
$13$ & $ \frac{1}{82899306 }\, \zeta(26) (46885126608855949+32769033282135949 m^2+23932490769301549 m^4 $ \\
 & $+16581373742311501 m^6+10657001539527466 m^8+6251551315383406 m^{10} $ \\
& $+3288926574739816 m^{12} +1518316248517036 m^{14}+597409832798656 m^{16}$ \\
& $+192183289614336 m^{18}  + 47355550989312 m^{20} +7940215990272 m^{22}+679111114752 m^{24})$\\ [.1cm]
$14$ & $ \frac{1}{6785560294 } \, \zeta(28) (18876870638218467389+13386158464071267389 m^2 $ \\
 & $+9936028183278303389 m^4+7037937839930870909 m^6+4662484330428335720 m^8$ \\
 & $+2849737068204025620 m^{10}+1584446206207320450 m^{12}+787977205890303150 m^{14} $ \\
 & $+343067669320476480 m^{16}+127034728547601280 m^{18}+38382891063494656 m^{20}$ \\
 & $+8866606540512256 m^{22}+1391299071164416 m^{24}+111174619856896 m^{26})$ \\ [.1cm]
$15$ & $ \frac{1}{6892673020804 }\, \zeta(30)(94347422281763617652599 +67773208432462419572599 m^2 $ \\
& $ +51020868251978707090999 m^4+36832275607229172708727 m^6 +25034877087339771576013 m^8$ \\
& $+15834856690259096336845 m^{10}+9212522702209442756325 m^{12} +4863979571266245490245 m^{14}$ \\
& $+2292234370375168188000 m^{16}+943898600252928069440 m^{18} +329986597479890425856 m^{20} $ \\
& $+93975868433656203776 m^{22}+20429529466374004736 m^{24} +3012284194329018368 m^{26} $ \\
& $ +225859109545705472 m^{28})$ \\ [.1 cm]
\hline
\end{tabular}
\vspace{0.3 cm}
\caption{The first 15 $(2m^2)^{v}S_{m,v}/ (m^2-1)$  for $v$, a positive integer}
\label{table3}\end{table}

As mentioned previously, values of the Gardner-Fisher sum have been
obtained in  Ref.\ \cite{K2011b} by an empirical approach. According to
this approach, the coefficients of $m^{-2i}$ in $S_{m,v}$ for $i<v$
were determined to be
\begin{equation}
C^{v}_{i}= c_{2v,i}\,\zeta(2v-2i) \Bigl( \frac{\pi}{2}\Bigr)^{2 i}\;,
\label{three-fourc}\end{equation}
while for $i=v$, they were given by
\begin{equation}
C^{v}_{v}= 2^{v} \Bigl(\frac{\pi}{2} \Bigr)^{2v} - \sum_{i=0}^{v-1}
2^{2v-2i} \,c_{2v,i} \Bigl(\frac{\pi}{2}\Bigr)^{2i}\,\zeta(2v-2i)\;.
\label{three-five}\end{equation}
In these results we have dropped dividing by $\zeta(2v)$, which occurs when this
factor is taken outside the results for $S_{m,v}$ to confirm Gardner's limit. In
addition, \eqref{three-five} has an extra factor of $2^v$ in the
first term on the rhs, which is missing in (44) of  Ref.\ \cite{K2011b}.
Moreover, in a future work it will be shown that the final
coefficient of $S_{m,v}$ can also be expressed as
\begin{equation}
C^{v}_{v}= - \frac{1}{2} \bigl( c_{2v,v} +1 \bigr) \Bigl( \frac{\pi}{2} \Bigr)^{2v} \;.
\label{three-six}\end{equation}
By equating like powers of $m^2$ between the above results and \eqref{three-three} we find that
\begin{equation}\label{three-seven}
c_{2v,i} = 2^{2i}\, \frac{ \Gamma(2v-2i)}{\Gamma(2v)} \, s(v,i)  \;,\quad i<v\;,
\end{equation}
\begin{equation} \label{three-sevena}
\sum_{i=0}^{v-1} \Bigl(2^{2v-2i}\, c_{2v,i} -\frac{2^{2i}}{(2v-1)!}\, s(v,i) \Gamma(2v-2i) \Bigr)
\Bigl( \frac{\pi}{2} \Bigr)^{2i-2v} \zeta(2v-2i) =2^{v} \;,
\end{equation}
and
\begin{equation}\label{three-sevenb}
\sum_{n=0}^{v-1} \pi^{2n}\,s(v,n) \,\Gamma(2v-2n) \, \zeta(2v-2n) =\frac{1}{2} \, \Gamma(2v) \bigl( c_{2v,v} +1 \bigr)
\Bigl( \frac{\pi}{2}\Bigr)^{2v} \;.
\end{equation}
These interesting results involving the generalized cosecant numbers have been verified by programming them in Mathematica.
Furthermore, by introducing \eqref{three-seven} into \eqref{kow2a}, we obtain the recurrence relation for the symmetric
polynomials, which is
\begin{equation*}
s(n+1,k+1) = s(n,k+1) + n^2 \, s(n,k) \;.
\end{equation*}

\section{Untwisted Dowker sum}
The integral approach of the previous section can be extended to the
situation where the external normalization factor in the
Gardner-Fisher sum is dropped and $\pi/2$ inside the trigonometric
power is replaced by $\pi/\ell$, where $\ell$ is any integer except
0. Then by using the same integral approach as before, we obtain
$$
\sum_{k=1}^{m-1}\csc^{2v}\left(\frac{k\pi}{\ell}\right)=\frac
1{\pi^{2v}\Gamma(2v)} \sum_{n=1}^{v-1}(2\pi)^{2n} s(v,n)
\sum_{k=1}^{m-1}\psi^{2v-2n-1}\left(\frac{2k}{\ell m}+\frac
{\ell-1}\ell\right).
$$
With the aid of No. 8.363(8) in  Ref.\ \cite{GR1994}, we can replace the
derivative of digamma function by the Hurwitz zeta function, thereby
arriving at
\begin{align} \label{eq4}
S_{m,v,\ell}=\sum_{k=1}^{m-1}\csc^{2v}\left(\frac{k\pi}{\ell}\right)&
=\frac{1}{\pi^{2v} \, \Gamma(2v)} \sum_{n=0}^{v-1}(2\pi)^{2n} s(v,n)
\sum_{k=1}^{m-1}\Gamma(2v-2n)
\nonumber\\
& \times \;\; \zeta^{2v-2n}\Bigl( 2v-2n,1+(2k/m-1)/|\ell| \Bigr).
\end{align}
For $|\ell| \neq 1$ or 2, the Hurwitz zeta function is intractable
and consequently, this generalization of the Gardner-Fisher sum will
not yield polynomials as in Table\ \ref{table3}. However, for the
important case of $\ell=1$, which is studied extensively in
Ref.\ \cite{Cv2007} and is known as the untwisted form of the Dowker sum
\cite{Do1992}, \eqref{eq4} reduces to
\begin{align}\label{eq5}
S_{m,v,1} = \sum_{k=1}^{m-1}\csc^{2v}\left(\frac{k\pi}{m}\right) &=
2^{2v+1} \sum_{n=0}^{v-1}\left(\frac{m}{2\pi}\right)^{2v-2n}
\frac{\Gamma(2v-2n)}{\Gamma(2v)}
\nonumber\\
& \times \;\; s(v,n)
\left(1-\frac{1}{m^{2v-2n}}\right)\zeta \bigl(2v-2n \bigl)\,.
\end{align}
Because there is no normalization factor outside the sum as in the
Gardner-Fisher sum, we obtain polynomials in powers of $m^2$ of
degree $v$. Although these polynomials are denoted by $C_{2v}(m)$ in
Ref.\ \cite{Cv2007}, we shall denote them by $q_{v}(m^2)$ with the
coefficients of $m^{2i}$ represented by $q_{v,i}$. From \eqref{eq5}
we find that
\begin{equation} \label{eq5a}
q_{v.0}=- 2^{2v+1} \sum_{n=0}^{v-1} (2\pi)^{2n-2v}\, \frac{\Gamma(2v-2n)}{\Gamma(2v)}\, s(v,n)\, \zeta(2v-2n)\;, \quad
q_{v,1}= \frac{1}{6} \, \frac{\Gamma(v) \,\Gamma(1/2)}{\Gamma(v+1/2)}\;,
\end{equation}
\begin{equation} \label{eq5b}
q_{v,i}= \frac{2^{2v-2i+1}}{\pi^{2i}}\, \frac{\Gamma(2i)}{\Gamma(2v)}\,s(v,v-i)\,\zeta(2i)\;, \;\;  i<v \,,
\quad {\rm and} \quad
q_{v,v}=\frac{2}{\pi^{2v}}\; \zeta(2v)\;.
\end{equation}
Moreover, we can use \eqref{three-seven} to express $q_{v,0}$ and $q_{v,i}$ in terms of the generalized cosecant numbers.
Hence, we obtain
$$q_{v,0}= -2 \sum_{n=0}^{v-1} \pi^{2n-2v}\, c_{2v,n} \,\zeta(2v-2n)\;,
$$
and
$$
q_{v,i}=\frac{2}{\pi^{2i}}\; c_{2v,v-i}\; \zeta(2i) \;, \quad i<v \, .
$$

\begin{table}
\small
\begin{tabular}{|c|l|} \hline
$v$ & $S_{m,v,1}/ (m^2-1)$ \\ [.1 cm]\hline
$1$ & $\frac{2\, \zeta(2)}{\pi^2}$ \\ [.1cm]
$2$  &  $\frac{2\, \zeta(4)}{\pi^4} (m^2+11)$  \\ [.1cm]
$3$ &  $\frac{\zeta(6)}{\pi^6} (2 m^4+23 m^2+ 191)$ \\ [.1cm]
$4$ & $ \frac{2\, \zeta(8)}{3 \pi^8} (m^2+11) (3m^4+10 m^2 +227)$ \\ [.1cm]
$5$ & $  \frac{\zeta(10)}{\pi^{10}} (2 m^8 +35 m^6+ 321 m^4 + 2125 m^2+ 14797)$ \\ [.1cm]
$6$ & $ \frac{2\,\zeta(12)}{1382 \, \pi^{12}} (1382 m^{10} + 28682 m^8 + 307961 m^6 + 2295661 m^4 + 13803157 m^2 + 92427157)  $ \\ [.1cm]
$7$ & $\frac{2\, \zeta(14)}{60 \pi^{14}} (60 m^{12} + 1442 m^{10} + 17822 m^8 + 151241 m^6 + 997801 m^4 + 5636617 m^2 + 36740617 )$\\ [.1cm]
$8$ & $ \frac{2\, \zeta(16)}{10851\, \pi{16}} (10851 m^{14} + 296451 m^{12} + 4149467 m^{10} + 39686267 m^8 + 292184513 m^6 $\\
& $ + 1777658113 m^4 + 9611679169 m^2 + 61430943169)$ \\ [.1cm]
$9$ & $ \frac{2\, \zeta(18)}{438670\, \pi^{18}} ( 438670 m^{16} + 13427317 m^{14} + 209998657 m^{12} + 2237483869 m^{10}+ 18276362179 m^8$ \\
& $ + 122248926511 m^6 + 701977130191 m^4 + 3674288164303 m^2 + 23133945892303)$ \\ [.1cm]
$10$ & $\frac{2\, \zeta(20)}{7333662 \pi^{20}}(7333662 m^{18} + 248602162 m^{16} + 4296730477 m^{14} + 50482177477 m^{12} $\\
 & $ + 453588822847 m^{10} + 3325534763347 m^8 + 20752800653227 m^6 + 114104919557227 m^4 $ \\
& $ + 582479437959787 m^2 + 3624331198599787)$ \\ [.1cm]
$11$ &  $\frac{2\, \zeta(22)}{4660980 \pi^{22}} (4660980 m^{20} + 173335206 m^{18} + 3280873486 m^{16} + 42140637451 m^{14} $ \\
& $ + 413195563231 m^{12} + 3299060965861 m^{10} + 22349066775541 m^8 + 132315729680101 m^6  $ \\
& $ + 703327910545381 m^4 + 3519581550481381 m^2 + 21691682977681381)$   \\ [.1cm]
$12$ & $\frac{2\, \zeta(24)}{7090922730 \pi^{24}}(7090922730 m^{22} + 287029381530 m^{20} + 5904893152686 m^{18} $\\
 & $+ 82314820117486 m^{16} + 874693313368831 m^{14} + 7557153052851631 m^{12} $ \\
& $ + 55300066644597091 m^{10} + 352631499074701891 m^8 + 2003157332714424931 m^6 $ \\
 & $ + 10363313134329413731 m^4 + 51026532637359173731 m^2+ 312017413700271173731) $\\ [.1cm]
$13$ & $ \frac{2\, \zeta(26)}{165798612 \pi^{26}} (165798612 m^{24} + 7256721342 m^{22} + 161222873682 m^{20} + 2424156289698 m^{18} $ \\
 & $+ 27750815238718 m^{16} + 257977419168313 m^{14} + 2028587745391093 m^{12}$ \\
& $+ 13879086707965453 m^{10} + 84366290294270413 m^8 + 463526111272448653 m^6  $ \\
& $+ 2345412070181900941 m^4 + 11394031603324326541 m^2 + 69213549869569446541 )$\\ [.1cm]
$14$ & $ \frac{2\, \zeta(28)}{ 13571120588 \pi^{28}}(13571120588 m^{26} + 638631887828 m^{24} + 15238841788898 m^{22} $ \\
 & $+ 245834814624698 m^{20} + 3016204736003030 m^{18} + 30020322332612430 m^{16}$ \\
 & $+ 252475090617525765 m^{14} + 1845413091251560365 m^{12} + 11967739987225201725 m^{10}$ \\
 & $+ 69975652378403124925 m^8 + 374033701594727629117 m^6 + 1857855957388613058877 m^4 $ \\
 & $+ 8923722772452603330877 m^2 + 53903636903066465730877)$ \\ [.1cm]
$15$ & $ \frac{2\, \zeta(30)}{13785346041608 \pi^{30}}( 13785346041608 m^{28} + 694064907756284 m^{26} + 17703093493737716 m^{24}  $ \\
& $ + 304993480084683806 m^{22} + 3992661121432093526 m^{20} + 42362161294746946250 m^{18}$ \\
& $+ 379446103825306975890 m^{16} + 2951191304716384278135 m^{14} +20345363828489173342455 m^{12} $ \\
& $+ 126302707637283630630775 m^{10} + 715104013050446845937527 m^8 $\\
& $ +3735238034142133535912311 m^6 + 18264356902365656743358839 m^4 $ \\
& $ + 86881942281626943067992439 m^2 + 522273861988577772410712439)$ \\ [.1 cm]
\hline
\end{tabular}
\vspace{0.3 cm}
\caption{The first 15 $S_{m,v,1}/ (m^2-1)$  for $v$, a positive integer}
\label{table4}\end{table}

As in the case of the Gardner-Fisher sum, the polynomials obtained from \eqref{eq5} possess a commnon factor of $(m^2\!-\!1)$. Consequently,
we can simplify the presentation of the polynomials for $S_{m,v,1}$ by removing this factor. Hence, Table\ \ref{table4} presents the
first 15 values of $S_{m,v,1}/(m^2-1)$, which were obtained by writing \eqref{eq5} as a one-line instruction in Mathematica as we did
for the Gardner-Fisher sum. In this instance the instruction becomes
\vspace{0.5cm}

CS[m$_{-}$, v$_{-}$] := (2$^{\wedge}$(2 v + 1)/(2 v - 1)!) Sum[(2 Pi/m)$^{\wedge}$(2 n -2 v)

SymmetricPolynomial[n, Table[k$^{\wedge}$ 2, \{k, 1, v - 1\}]]

Gamma[2 (v - n)] Zeta[2 (v - n)] (1 - 1/(m$^{\wedge}$(2 (v - n)))), \{n, 0, v - 1\}]
\vspace{0.5cm}
\newline
Then we can use the same intructions below the instruction for the
Gardner-Fisher sum except that S[m,v] is now replaced by CS[m,v] to
tabulate the polynomials and time the calculation. The results in
Table\ \ref{table4} took 0.12 seconds to compute as the same Venom
laptop mentioned previously. It should also be mentioned that the
first five results in the table are identical to those given in
Ref.\ \cite{Cv2007}. In fact, these authors prove that
\begin{equation}
q_v \left( m^2 \right)=  (-1)^{v-1} \; \frac{2^{2n}}{(2n)!} \sum_{n=0}^v
\binom{2v}{2n} \,B_{2v-2n}\, B^{(2v)}_{2n}(v)\, m^{2v-2n}\;,
\label{three-eight}\end{equation} where $B_{2v-2n}$ are the ordinary
Bernoulli numbers and $B^{(m)}_{k}(x)$ are the Bernoulli polynomials
of order $m$ and degree $k$ and sometimes referred as to N\"{o}rlund
polynomals. By equating like powers of $m$ between \eqref{eq5} and
\eqref{three-eight}, we obtain the following results:
\begin{equation}
B^{(2v)}_{2n}(v)= (-1)^n \, (2n)! \,\frac{\Gamma(2v-2n)}{\Gamma(2v)}\,s(v,n) \,, \quad n<v\;,
\label{three-nine}\end{equation}
and
\begin{equation}
B^{(2v)}_{2v}(v)= (-1)^v \,4 v \sum_{n=0}^{v-1} (2\pi)^{2n-2v}\,
\Gamma(2v-2n) \, s(v,n)\, \zeta(2v-2n) \;.
\label{three-ten}\end{equation} In obtaining this result we have
used No.\ 9.616 in  Ref.\ \cite{GR1994}, which expresses the Bernoulli
numbers in terms of the Riemann zeta function. Furthermore,
introducing \eqref{three-seven} into \eqref{three-nine} and
\eqref{three-ten} yields
$$
B^{(2v)}_{2n}(v)= (-1)^n \,2^{-2n} \,(2n)! \,c_{2v,n} \,, \quad n<v \;,
$$
and
$$
B^{(2v)}_{2v}(v)= (-1)^v \, 2^{1-2v}\, \Gamma(2v+1) \sum_{n=0}^{v-1} \pi^{2n-2v} \,c_{2v,n}\, \zeta(2v-2n)\;.
$$
The above results can be verified in Mathematica \cite{W1992}, where
the Bernoulli polynomials of order $m$ and degree $k$ are determined
by using the NorlundB routine.

\section{Other Sums}
We can use the results of the previous sections to consider more intricate
trigonometric power sums than the Gardner-Fisher and untwisted Dowker sums.
E.g., consider the following series:
\begin{equation}
S^{CC}_{m,v,w,l} := \sum_{k=1}^{m-1} \cot^{2v}\Bigl( \frac{k \pi}{\ell m} \Bigr)
\csc^{2w} \Bigl( \frac{k \pi}{\ell m} \Bigr)\;\;,
\label{four-one}\end{equation}
where $v \geq 0$ and $w \geq 0$ excluding $v+w \!=\! 0$, while $\ell \!=\! 1$ and
$\ell \!=\! 2$ correspond respectively to the untwisted Dowker and Gardner-Fisher cases.
We could introduce a factor of $\cos \left(2 a k \pi/\ell m\right)$, where $a$ is an
integer less than $\ell m-1$, into the summand. For $\ell=1$, this becomes an
extension of the twisted Dowker sum \cite{Do1992}, which we aim to study in a future
work.

The above trigonometric power sum can also be expressed as
\begin{equation}
S^{CC}_{m,v,w,\ell}= \sum_{k=1}^{m-1} \frac{\cos^{2v}(k \pi/\ell m)}{\sin^{2v+2w}(k \pi/\ell m)}\;.
\label{four-two}\end{equation}
Now we replace the cosine power in the numerator by $\left(1-\sin^2(k\pi/\ell m)\right)^v$ and apply the
binomial theorem \cite{wi2015}, thereby obtaining
\begin{equation}
 S^{CC}_{m,v,w,\ell}= \sum_{j=0}^{v} \sum_{k=1}^{m-1} (-1)^{v-j} \, \binom{v}{j} \,\csc^{2w+2j}
\Bigl( \frac{k \pi}{\ell m} \Bigr) \;\;.
\label{four-three}\end{equation}
Therefore, the sum represents a finite sum of Fisher-Gardner and untwisted Dowker sums
depending upon the value of $\ell$.

If $w+v \leq 15$, then we can use the results displayed in Tables\ \ref{table3} and \ref{table4} to
obtain the values of $S^{CC}_{m,v,w,\ell}$. Moreover, denoting the quantities listed in Table\
\ref{table3} as $R_v(m^2)$, i.e.\ $R_v \left( m^2 \right)= (2 m^2)^v S_{m,v}/(m^2-1)$, we find
that $S^{CC}_{m,v,w,2}$ reduces to
\begin{equation}
 S^{CC}_{m,v,w,2}= \left(m^2 -1\right) \sum_{j=0}^{v}(-1)^{v-j} \binom{v}{j}
\frac{2^{w+j}}{\pi^{2w+2j}}\, R_{w+j}\left(m^2\right) \;\;.
\label{four-four}\end{equation}
If $v$ and $w$ are set equal to 5 and 4 respectively, then using the values of $R_4 \left(m^2 \right)$ to
$R_9\left(m^2 \right)$ in Table\ \ref{table3} and the above result, we find that
\begin{align}
S^{CC}_{m,5,4,2} & = \frac{16}{194896477400625} \, \left(m^2-1 \right) \left( 4m^2-1 \right) \left( 2280413161 + 712556555 m^2
\right.
\nonumber\\
 & - \;\; 2906805048 m^4 - 2535353600 m^6 + 2920623488 m^8 + 2565749760 m^{10}
\nonumber\\
& \left. - \;\;3310462976 m^{12} +  898396160 m^{14} \right)  \;\;.
\label{four-five}\end{align}
A general formula for the $\ell=2$ case can be obtained by using Eqs.\ (\ref{three-fourc}) to (\ref{three-six}),
which represent the coefficients of $S_{m,v}$. Then we arrive at
\begin{equation}
S^{CC}_{m,v,w,2}=  \sum_{j=0}^{v} (-1)^{v-j} \binom{v}{j} \, \frac{2^{2v+2j}}{\pi^{2v+2j}} \, \sum_{i=0}^{v+j}
C^{v+j}_i m^{2v+2j-2i} \;\;.
\label{four-six}\end{equation}
The above result can also be expressed in terms of the symmetric polynomials of quadratic integers since from
\eqref{three-fourc}, \eqref{three-six}, \eqref{three-seven} and \eqref{three-sevenb}, we have
\begin{equation*}
C^v_i= (2v)_{-2n}\, \zeta(2v-2i)\, s(v,i)\;, \quad i<v\;,
\end{equation*}
and
\begin{equation*}
C^v_v=-\sum_{n=0}^{v-1} \pi^{2n} (2v)_{-2n}\, \zeta(2v-2n)\, s(v,n) \;.
\end{equation*}

For $\ell=1$ or the untwisted Dowker case, we can use the results of Table\ \ref{table4} to obtain
the values of the sums, $S^{CC}_{m,v,w,1}$, when $v+w \leq 15$. If we denote the quantities listed in the table
as $T_v \left( m^2 \right)$, i.e.\ $T_v\left( m^2 \right) =S_{m,v,1}/(m^2-1)$, then we find that
\begin{equation}
S^{CC}_{m,v,w,1}= \left(m^2 -1\right) \sum_{j=0}^{v}(-1)^{v-j} \binom{v}{j}
T_{w+j}\left( m^2 \right) \;.
 \label{four-seven}\end{equation}
For $v=6$ and $w=3$, the above result with the aid of the respective values in Table\ \ref{table4} yields
\begin{align}
\sum_{k=1}^{m-1} \cot^{12}\Bigl( \frac{k \pi}{ m} \Bigr) & \csc^{6} \Bigl( \frac{k \pi}{ m} \Bigr) =
\frac{1}{194896477400625} \left(m^2-1 \right) \left( m^2-4 \right)
\left( -29787342748 \right.
\nonumber\\
& + \;\; 1960688815 m^2 +  3595494399 m^4 - 275848135 m^6 -  335395979 m^{8}
\nonumber\\
& \left.  + \;\; 98107275 m^{10} - 10795297 m^{12} + 438670 m^{14} \right) \;.
\label{four-eight}\end{align}
Both \eqref{four-five} and \eqref{four-eight} have been checked for specific values of $m$ by
calculating the decimal values of the sums on the lhs and comparing them with the decimal values
of the rational quantities on the rhs's in Mathematica. In terms of the $q_{v,i}$ given by
\eqref{eq5a} and \eqref{eq5b}, we can express $S^{CC}_{m,v,w,1}$ more generally as
\begin{align}
S^{CC}_{m,v,w,1}= \sum_{j=0}^{v}(-1)^{v-j} \binom{v}{j} \sum_{i=0}^{w+j}q_{w=j,i} \, m^{2i} \;.
\label{four-nine}\end{align}

Another trigonometric power sum that can be studied with the aid of the previous sections is:
\begin{equation}
S^{TS}_{m,v,w,l} := \sum_{k=1}^{m-1} \tan^{2v}\Bigl( \frac{k \pi}{\ell m} \Bigr)
\sec^{2w} \Bigl( \frac{k \pi}{\ell m} \Bigr)\;\;,
\label{four-ten}\end{equation}
where the same conditions apply on $v$ and $w$ as in \eqref{four-one}. For $\ell=2$, if we replace
$k$ by $m-k$, then we find that $S^{TS}_{m,v,w,2}=S^{CC}_{m,v,w,2}$. Hence, there is no need to consider
this case.

On the other hand, for $\ell \!=\! 1$, there is a possibility that one of the values of $k$ can produce powers
of $\cos(\pi/2)$ in the denominator when $m$ is an even integer. Therefore, this value of $k$ needs to be
excluded when $m=2n$, where $n$ is a non-zero integer. Therefore, we modify the above sum to
\begin{equation}
S^{TS}_{2n,v,w,1} = \sum_{k=1 \atop{k \neq n}}^{2n-1} \tan^{2v}\Bigl( \frac{k \pi}{2n } \Bigr)
\sec^{2w} \Bigl( \frac{k \pi}{2n} \Bigr)\;\;,
\label{four-eleven}\end{equation}
We now split the above sum into two separate sums, the first ranging from $k=1$ to $n-1$ and the second,
from $k=n-1$ to $2n$. In the second sum we replace $k$ by $2n-k$, which yields the first sum again.
Consequently, \eqref{four-eleven} reduces to
\begin{equation}
S^{TS}_{2n,v,w,1} = 2 \sum_{k=1}^{n-1} \tan^{2v}\Bigl( \frac{k \pi}{2n } \Bigr)
\sec^{2w} \Bigl( \frac{k \pi}{2n} \Bigr)\;\;,
\label{four-twelve}\end{equation}
By replacing $k$ by $n \!-\! k$, we find that $S^{TS}_{2n,v,w,1}=2 S^{CC}_{n,v,w,2}$. The case of $m \!=\! 2n \!+\! 1$
can also be reduced, but it yields a second trigonometric power sum with an alternating summand. Such sums will be studied
in a future work. Finally, we add that the $\ell \!=\! 1$, $w \!=\! 0$ and $m \!=\! 2n \!+\!1$ case of \eqref{four-ten}
has been studied by Shevelev and Moses in Ref.\ \cite{Sh2014}, where they give the polynomial values of the sum
for the first five values of $v$.

\section{Conclusion}
In this paper we have presented a new integral approach for
evaluating the Gardner-Fisher sum or $S_{m,v}$ as defined by either
\eqref{one} or \eqref{onea}, which is not only computationally
expedient compared with other methods, but has also enabled us to
quantify the coefficients of the polynomials as evidenced by
\eqref{three-four} to \eqref{three-fourb}. In so doing, we have been
able to relate the generalized cosecant numbers of Ref.\ \cite{K2011a} to
the symmetric polynomials $s(v,n)$ over the set of quadratic powers,
$\{1,2^2,3^2,\dots,(v-1)^2\}$, via \eqref{three-seven}, which was
achieved by matching the general result given here by
\eqref{three-three} with the empirically determined results of (43)
and (44) in Ref.\ \cite{K2011b}.

To demonstrate the versatility of our integral approach, it was then
extended to situations where $\pi/2$ in the trigonometric power of
the Gardner-Fisher sum was replaced by $\pi/\ell$. As a consequence,
we were able to evaluate the sums for the $\ell \!=\! 1$ case or
$S_{m,v,1}$, which is known as the untwisted Dowker sum
\cite{Do1992} and has been studied extensively by Cvijovi\'c and
Srivastava \cite{Cv2007}. The latter authors obtain a general result
for the sum, which is given by another unwieldy sum whose summand is
a product of Bernoulli numbers and the esoteric N\"{o}rlund
polynomials as in \eqref{three-eight}. As a result, one is unable to
ascertain the mathematical forms for the coefficients of the
polynomials in $S_{m,v,1}$. Hence, one has to rely on a software
package such as Mathematica to generate the final forms for the
untwisted Dowker sum. Nevertheless, the five values presented in
\cite{Cv2007} agree with the first five results in Table\
\ref{table4}. Furthermore, by comparing their form for the untwisted
Dowker sum with our \eqref{eq5}, we are able to express the
particular values of the N\"{o}rlund polynomials in their result
either in terms of the specific symmetric polynomials $s(v,n)$
presented here or in terms of the generalized cosecant numbers.
By using the results of the previous sections we were able to sketch out
the calculations for more intricate sums involving products of powers of
cotangent and tangent with powers of cosecant and secant respectively.

This paper, which has resulted in a cross-fertilization of the fields of
classical analysis, number theory and computational/experimental mathematics,
represents the introductory work of a more ambitious
programme where the ideas and methods presented here and in
\cite{K2011b,K2011a} are to be extended beyond the trigonometric power
sums appearing in Refs.\ \cite{Be2002} and \cite{Cv2012}.
Included in this investigation will be the cases where summands
alternate in sign. Once again, the existing results give similar unwieldy
results to \eqref{three-eight} and thus, do not provide the
interesting mathematics associated with the coefficients of the
resultant polynomials.

\section{Acknowledgment}
The authors are very grateful to the Office of the Vice-President for Research, Kuwait University for supporting and funding
this project under Kuwait University Research Grant No. SM03/13.

\end{document}